\documentclass[12pt]{article}

\usepackage{amsfonts}  

\title{Eigenvalue and Dirichlet problem for fully-nonlinear operators in non smooth domains.}

\author{I. Birindelli, F. Demengel}

\date{}

\catcode`@=11 \@addtoreset{equation}{section} \catcode`@=12

\newtheorem{theo}{Theorem}[section]
\newtheorem{prop}[theo]{Proposition}
\newtheorem{rema}[theo]{Remark}
\newtheorem{defi}[theo]{Definition}
\newtheorem{cor}[theo]{Corollary}

\def\R{{\rm I}\!{\rm  R}}
\def\grad{\nabla}

\begin{document}
\maketitle

\section{Introduction}
In this paper we study the maximum principle, the existence of eigenvalue and the existence of solution for the  Dirichlet problem for operators which are fully-nonlinear, elliptic but presenting some singularity or degeneracy which are similar to those of the  p-Laplacian, the novelty resides in the fact that we consider the equations in bounded domains which only satisfy the exterior cone condition.

Before defining the precise notions described above let us recall that Berestycki, Nirenberg and Varadhan in \cite{BNV}, have proved maximum principle, principal eigenvalue and Dirichlet problem for linear uniformly elliptic operators $Lu={\rm tr} A(x)D^2u+b(x)\cdot\grad u+c(x)u$ in domains without any regularity condition on the boundary.

In order to do so, they need to define the concept of boundary condition. Hence, using Alexandrov Bakelman Pucci inequality and Krylov-Safonov Harnack's inequality  they first prove the existence of $u_o$ a strong solution  of 
$${\rm tr} A(x)D^2u_o+b(x)\cdot\grad u_o=-1\quad \mbox{in}\quad \Omega$$
which is zero on the points of the boundary that have some smoothness. Then they define the boundary condition  for the full operator $L$ through this function $u_o$.
Their paper, which constructs the principal eigenvalue using only the maximum principle, has allowed to 
generalize the notion of eigenvalue to fully-nonlinear operators, see e.g. \cite{BEQ,QS,IY,J,BD2,BD3,P}.

Here, as in \cite{BD2,BD3,P} we shall consider operators that satisfy:

(H1)
 $F: \Omega\times \R^N\setminus\{0\}\times S\rightarrow\R$, 
and  $\forall t\in \R^\star$, $\mu\geq 0$,
 $F(x, tp,\mu X)=|t|^{\alpha}\mu F(x, p,X)$.

\bigskip

(H2)
There exist $0< a<A$, 
for $x\in \overline{\Omega}$, $p\in \R^N\backslash \{0\}$, $M\in S$,  $N\in S$, 
$N\geq 0$
$$
a|p|^\alpha tr(N)\leq F(x,p,M+N)-F(x,p,M) \leq A
|p|^\alpha tr(N).
$$
and other "regularity" conditions.

For this class of operators it is not known whether the Alexandrov Bakelman Pucci inequality holds true, hence in our previous works  \cite{BD3, BD2} we supposed that  $\partial\Omega$ was ${\cal C}^2$.  The regularity of the boundary in those papers played a crucial role because it allowed to use the distance function to construct sub and super solutions. This was the key step in the proof of the maximum principle.
Here, instead, we shall suppose that $\Omega$ satisfies only the "uniform exterior" cone condition i.e.

\noindent {\em There exist $\psi>0$ and $\bar r>0$ such that    for any $z\in\partial \Omega$ and  for an axe through $z$ of direction $\vec n$, }
$$ Co:=\{x:\ \frac{(x-z)\cdot\vec n}{|z-x|}\leq \cos\psi\}, \quad\quad  Co\cap\overline{\Omega}\cap B_{\bar r}(z)=\{z\}.$$
This cone condition allows to construct some barriers and consequently a function which will play the same role as $u_o$ in \cite{BNV}.
In particular we can prove that there exists an eigenfunction $\varphi>0$, solution of
$$\left\{ \begin{array}{cc}
F(x,\grad\varphi, D^2\varphi)+h(x)\cdot\nabla \varphi|\nabla \varphi|^\alpha+  (V(x)+\bar \lambda(\Omega)) \varphi^{1+\alpha}=0 & {\rm in }\ \Omega\\
 \varphi = 0  & {\rm on}\ \partial \Omega
 \end{array}\right.$$
 for 
\begin{eqnarray*}
\bar\lambda(\Omega)&=& \sup \{ \lambda, \exists \ u>0\ {\rm in} \ \Omega ,\\
&& F(x,\grad u,D^2u)+h(x)\cdot\nabla u|\nabla u|^\alpha+  (V(x)+ \lambda) u^{1+\alpha} \leq 0\   {\rm in} \ \Omega \}.
\end{eqnarray*}
Finally  in the last section we also  define
 $$\lambda_e = \sup\{\lambda (\Omega^\prime), \ \Omega\subset \subset \Omega^\prime,\ \Omega^\prime \ {\rm regular\ and \ bounded}\}$$
and 
\begin{eqnarray*}
\tilde \lambda&=& \sup \{ \lambda, \exists \ u>0\ {\rm in} \ \overline{\Omega} , \\
& &F(x,\grad u,D^2u)+h(x)\cdot\nabla u|\nabla u|^\alpha+  (V(x)+ \lambda) u^{1+\alpha}\leq 0\}.
\end{eqnarray*}

We  prove that 
$\lambda_e = \tilde \lambda$ 
and  that this value  is an "eigenvalue" in the sense that there exists some $\phi_e >0$, which satisfies 
$$
\left\{\begin{array}{lc}
F(x,\nabla \phi_e, D^2\phi_e) + h(x)\cdot \nabla \phi_e |\nabla \phi_e|^\alpha + (V(x)+ \lambda_e(\Omega))\phi_e^{1+\alpha}=0 & \mbox{in}\quad\Omega\\
\phi_e=0 & \mbox{on}\quad\partial\Omega.
\end{array} 
\right.$$
We also prove that for any $\lambda<\lambda_e$ the maximum principle holds and there exists a solution of  the Dirichlet problem when the right hand side is negative. 

Observe that $\lambda_e  \leq \overline{ \lambda}$,  and furthermore if $\Omega$ is smooth, the equality holds.   It is an open problem to know if if the equality still holds when $\Omega$ satisfies only the exterior cone condition.
Let us observe that the  identity of these values is equivalent to the existence of a maximum principle for $\lambda < \bar\lambda$.

\section{Assumptions on $F$}

The following hypothesis will be considered

\begin{itemize}

\item[(H1)] 
 $F: \Omega\times \R^N\setminus\{0\}\times S\rightarrow\R$, 
and  $\forall t\in \R^\star$, $\mu\geq 0$,
 $F(x, tp,\mu X)=|t|^{\alpha}\mu F(x, p,X)$.

\item[(H2)]

There exist $0< a<A$,
for $x\in \overline{\Omega}$, $p\in \R^N\backslash \{0\}$, $M\in S$,  $N\in S$, 
$N\geq 0$

\begin{equation}\label{eqaA}
a|p|^\alpha tr(N)\leq F(x,p,M+N)-F(x,p,M) \leq A
|p|^\alpha tr(N).
\end{equation}

\item[(H3)]
There exists a continuous function $\tilde \omega$, $\tilde \omega(0)=0$
such that for all $(x,y)\in \Omega^2 $,  $\forall p\neq 0$, $\forall X\in S$
$$|F(x,p,X)-F(y,p,X)|\leq \tilde \omega(|x-y|) |p|^\alpha |X|.$$
\item [(H4)]
 There exists a
continuous function $  \omega$ with $\omega (0) = 0$, such that if
$(X,Y)\in S^2$ and 
$\zeta\in \R^+$ satisfy
$$-\zeta \left(\begin{array}{cc} I&0\\
0&I
\end{array}
\right)\leq \left(\begin{array}{cc}
X&0\\
0&Y
\end{array}\right)\leq 4\zeta \left( \begin{array}{cc}
I&-I\\
-I&I\end{array}\right)$$
and $I$ is the identity matrix in $\R^N$,
then for all  $(x,y)\in \R^N$, $x\neq y$
$$F(x, \zeta(x-y), X)-F(y,  \zeta(x-y), -Y)\leq \omega
(\zeta|x-y|^2).$$
\end{itemize}

\begin{rema} When no ambiguity arises we shall sometime write 
$F[u]$ to signify $F(x,\grad u,D^2u)$.
\end{rema}

We assume that $h$  and $V$ are some continuous and bounded functions on $\overline\Omega$ and 

(H5)
- Either $\alpha \leq 0$ and $h$  is H\"older continuous  of exponent $1+\alpha$,

\hspace{1cm} - or  $\alpha >0$ and 
$h(x)-h(y)\cdot x-y\leq 0$

\bigskip

The solutions that we consider will be taken in the sense of viscosity.  For convenience of the reader 
we state the precise definition.

\begin{defi}\label{def1}

 Let $\Omega$ be a bounded domain in
$\R^N$, then
$v$,   bounded and continuous on $\overline{\Omega}$ is called a viscosity super solution (respectively sub-solution)
of
$F(x,\grad u,D^2u)+h(x)\cdot\grad u|\grad u|^\alpha=f(x,u)$ if for all $x_0\in \Omega$, 

-Either there exists an open ball $B(x_0,\delta)$, $\delta>0$  in $\Omega$
on which 
$v= cte= c
$ and 
$0\leq f(x,c)$, for all $x\in B(x_0,\delta)$ (respectively $0\geq f(x,c)$)

-Or
 $\forall \varphi\in {\mathcal C}^2(\Omega)$, such that
$v-\varphi$ has a local minimum on $x_0$ (respectively a local maximum ) and $\grad\varphi(x_0)\neq
0$, one has
$$
F( x_0,\grad\varphi(x_0),
 D^2\varphi(x_0))+h(x_0)\cdot\grad \varphi(x_0)|\grad \varphi(x_0)|^\alpha\leq f(x_0,v(x_0)).
$$
(respectively
$$
F( x_0,\grad\varphi(x_0),
 D^2\varphi(x_0))+h(x_0)\cdot\grad \varphi(x_0)|\grad \varphi(x_0)|^\alpha\geq f(x_0,v(x_0)).)
 $$

\end{defi}

We now recall what we mean by first eigenvalue and some of  the properties  of this eigenvalue. 

When $\Omega$ is a bounded domain we define 
$$\bar\lambda (\Omega)= \sup \{ \lambda, \exists \ \varphi>0\ {\rm in} \ \Omega , F[\varphi]+h(x)\cdot\nabla \varphi|\nabla \varphi|^\alpha+  (V(x)+ \lambda) \varphi^{1+\alpha} \leq 0\}$$

When $\Omega$ is a bounded regular set, we proved in \cite{BD1} that 

\begin{theo}\label{exphi}

Suppose that $F$ satisfies (H1)-(H4), that $h$ satisfies (H5), and that $V$ is continuous and bounded. 
Suppose that $\Omega$ is a bounded regular domain. 

There exists $\varphi>0$ which is a solution of 
$$\left\{ \begin{array}{cc}
F[\varphi]+h(x)\cdot\nabla \varphi|\nabla \varphi|^\alpha+  (V(x)+\bar \lambda) \varphi^{1+\alpha}=0 & {\rm in }\ \Omega\\
 \varphi = 0  & {\rm on}\ \partial \Omega.
 \end{array}\right.$$
Moreover $\varphi$ is strictly positive inside $\Omega$ and is H\"older continuous \end{theo}
We now recall some  properties of the eigenvalue : 

\begin{theo}\label{pmaxreg}

 Suppose that $\Omega$ is a bounded regular domain , and that $F$,  $h$,  and $V$ satisfy the previous assumptions. 
  Suppose that $\lambda < \bar\lambda$ and that $u$ satisfies 
  $$\left\{ \begin{array}{lc}
  F(x, \nabla u, D^2 u) + h(x) \cdot \nabla u  |\nabla u|^\alpha + (V(x)+\lambda) |u|^\alpha u \geq 0&{\rm in } \ \Omega\\
  u\leq 0 & \ {\rm on } \ \partial \Omega.
  \end{array}\right.$$
  Then $u\leq 0$ in $\Omega$.
  \end{theo}

 We now recall the following comparison principle which holds without assumptions on the regularity of the bounded domain $\Omega$, 

\begin{prop}\label{comparisonprinciple}
Suppose that $\beta(x,.)$ is non decreasing and $\beta(x,0)=0$, that $w$ is an upper semicontinuous sub-solution of 
$$F(x, \nabla w, D^2 w) + h(x)\cdot \nabla w |\nabla w|^\alpha-\beta (x, w(x)) \geq g$$
and  $u$ is a lower semicontinuous  supersolution of 
$$F(x,\nabla u, D^2 u)+ h(x)\cdot \nabla u |\nabla u|^\alpha-\beta (x, u(x)) \leq f$$
with $g$ lower  semicontinuous, $f$ upper semicontinuous ,  $f<g$ in $\Omega $
and 
$$\limsup ( w(x_j)-u(x_j))\leq 0,$$
for all $x_j\rightarrow \partial \Omega$. 
Then 
$w\leq u$ in $\Omega$.
\end{prop}

\begin{rema}
The result still holds if $\beta$ is increasing and $f\leq g$ in $\Omega$.
\end{rema}

The proof is as in \cite{BD1}.

We also recall the following weak  comparison principle.
 
 \begin{theo}\label{thcomp}
Suppose that $\Omega$ is some bounded open set. Suppose that $F$ satisfies (H1), (H2), and (H4),  
that $h$ satisfies (H5) and $V$ is continuous and
bounded.  Suppose that   $f\leq 0$, $f$ is
upper semi-continuous and
$g$ is lower semi-continuous  with $f\leq  g$.

Suppose that there exist
$u$ and
$v$  continuous ,  $v\geq 0$, satisfying
\begin{eqnarray*}
 F(x, \nabla u, D^2 u)+ h(x)\cdot \nabla u |\nabla u|^\alpha +V(x) |u|^\alpha u& \geq & g\quad 
\mbox{in}\quad \Omega \\ 
F(x,  \grad v,D^2 v)
+ h(x)\cdot \nabla v |\nabla v|^\alpha+ V(x)) v^{1+\alpha} & \leq & f \quad  \mbox{in}\quad
\Omega  \\ 
u \leq  v &&   \quad  \mbox{on}\quad \partial\Omega.
\end{eqnarray*}
Then $u \leq v$ in $\Omega$ in each of these two cases:

\noindent 1) If $v>0$ on $\overline{ \Omega}$ and either $f<0$ in 
$\Omega$, 
 or  $g(\bar x)>0$ on every point $\bar x$ such that  $f(\bar x)=0$,

\noindent  2) If $v>0$ in $\Omega$, $f<0$ in $\overline{\Omega}$  and $f<g $ on
$\overline\Omega$.
\end{theo}

\section{Barriers in non smooth domains}
 
 In this section we assume that $\Omega$  satisfies the exterior cone condition.  More precisely we assume that there exists $\bar r$ and $\psi\in ]0,\pi[$ such that for each $z\in\partial\Omega$  the set $\Omega\cap  B(z, \bar r)$ is included in the open cone 
which, up to change of coordinates can be given by
$$T_\psi= \{ r\in ]0, \bar r[, 0\leq \arccos\left( {x_N\over r} \right)\leq \psi\}$$
choosing the main direction of the cone to be $e_N$.  Indeed, in that case, the exterior of $\Omega$  contains at least the set of $(x^\prime, x_N)$ with $-1\leq {x_N\over r} \leq \cos \psi$, $r< \bar r$. 

On the operator $F$ we suppose that it satisfies conditions (H1), (H2), (H3) and (H4), while $h$ satisfies
(H5).
\subsection{Local  barriers}

Under the exterior cone condition we are going to construct a local barrier i.e. for any $z\in\partial\Omega$,   
a supersolution   in a neighborhood of $z$, of $F[v]+ h(x)\cdot \nabla v|\nabla v|^\alpha \leq -1$,    such that $c|x-z|^\gamma \leq v(x)\leq C |x-z|^\gamma$ for some $\gamma\in(0,1] $ and for some constant $c$ and $C$  which depend on $\psi, a, A,\gamma,\bar r$. This barrier is constructed on  the model of those given by Miller for the Pucci operators in \cite{M1}.

 We define   
$$v = |x-z|^\gamma \varphi(\theta)$$
where 
$\theta = \arccos \left({x_N-z_N\over |x-z|}\right)$. Without loss of generality, we suppose that $z=0$.

We suppose first that $h\equiv 0$ and at the end of the proof, we will say which are  the changes that need to be done when  $h\not\equiv 0$ . 
We shall first  show that there exists  $\varphi$ a solution of some differential linear equation such that   $v$ is a super solution of
$$F(x,\nabla v , D^2 v)\leq -b$$
where $b$ is a positive constant that depends only on $\psi$, $\gamma$,  $r_o$ and the structural constant of the operator.
 It will be useful for the following to observe that $1\geq	\frac{x_N}{r}\geq \cos \psi$ on the considered set.

  Let $x= (x_1, \cdots, x_N)=(x',x_N)$. Let $r=|x|$ and $r'=|x'|$.
  We shall also use the following notation $X'=(x',0)$. 
  
One has:
$$\nabla v=\gamma r^{\gamma-2}x\varphi(\theta)+ r^\gamma \varphi^\prime(\theta)\grad \theta$$
and
\begin{eqnarray*}
D^2v&=& r^{\gamma-2}\varphi\gamma\left(I+\frac{(\gamma-2)}{r^2}x\otimes x\right)\\
&&+r^{\gamma-2}\varphi^\prime\left(r^2D^2\theta +\gamma (\grad\theta\otimes x+x\otimes\grad\theta)\right)\\
&&+r^{\gamma-2}\varphi^{\prime\prime}\left(r^2\grad\theta\otimes\grad\theta\right).
\end{eqnarray*}
We now suppose that $\varphi\geq 0$, $\varphi^\prime\leq 0$ and $\varphi^{\prime\prime}\leq 0$ then

\begin{eqnarray*}
{\cal M}^+_{a,A}(D^2v)&\leq& r^{\gamma-2}\left(\varphi\gamma{\cal M}^+_{a,A}\left(I+\frac{(\gamma-2)}{r^2}x\otimes x\right)+\right.\\
&&+\ \varphi^\prime{\cal M}^-_{a,A}(r^2D^2\theta +\gamma (\grad\theta\otimes x+x\otimes\grad\theta))\\
&&\left.+\ \varphi^{\prime\prime}{\cal M}^-_{a,A}(r^2\grad\theta\otimes\grad\theta)\right).
\end{eqnarray*}
Since we need to find the eigenvalues of the above matrices let us remark that
 $$\nabla \theta = \frac{1}{r'}(\frac{x_Nx}{r^2}-e_N)={x^\perp\over r^2}$$
 with
 $$x^\perp =\frac{ x_N}{r'}X'-r'e_N=\cot \theta x-{r^2\over r^\prime} e_N.$$
 In particular $x^\perp\cdot x=0$ and $|x^\perp|=r$ and we obtain
 
 \begin{eqnarray*}
 {\cal M}^+_{a,A}\left(I+\frac{(\gamma-2)}{r^2}x\otimes x\right)&=&A(N-1)+a(\gamma-1)\\
 {\cal M}^-_{a,A}(r^2\grad\theta\otimes\grad\theta)&=&ar^2|\grad\theta|^2=a\\
 {\cal M}^-_{a,A}(\gamma (\grad\theta\otimes x+x\otimes\grad\theta))&=&\gamma|\grad\theta|r(a-A)=\gamma(a-A).
 \end{eqnarray*}
 To complete the calculation we need to  compute
 \begin{eqnarray*}
D^2\theta&=&-\frac{1}{r'^2}\frac{X'}{r'}\otimes\left(\frac{x_Nx}{r^2}-e_N\right )+\\
&&+\frac{1}{r'}\left[\frac{1}{r^2}e_N\otimes x-\frac{2x_N}{r^4}x\otimes x +\frac{x_N}{r^2}I\right]\\
&=&-\frac{x_N}{(r')^3 r^2}X'\otimes X' +\\
&& +\frac{1}{r'r^2}\left[X'\otimes e_N+e_N\otimes X' +x_N e_N\otimes e_N -2\frac{x_N}{r^2} x\otimes x +x_NI\right].
\end{eqnarray*}
To estimates the eigenvalues  of $r^2D^2\theta$ we shall use the following facts and notations:

$I_{N-1}$ indicate the identity $(N-1)\times(N-1)$ matrix,

$$ I'=\left(\begin{array}{cc}
I_{N-1} &0\\
0  & 0
\end{array}\right), 
\quad I=I'+e_N\otimes e_N,$$
$$\quad x\otimes x =X'\otimes X'+x_N(X'\otimes e_N+e_N\otimes X' )+x_N^2e_N\otimes e_N.
$$
Then
 \begin{eqnarray*}
r^2D^2\theta&=& \frac{x_N}{r'}\left(-\frac{1}{(r')^2}X'\otimes X' +I'\right)+ \frac{x_N}{r'}(2-2\frac{x_N^2}{r^2})e_N\otimes e_N\\
&&+\ \frac{1}{r'}(1-2\frac{x_N^2}{r^2})(X'\otimes e_N+e_N\otimes X' )-2\frac{x_N}{r'}\frac{1}{r^2}X'\otimes X'.
\end{eqnarray*}
One has 
\begin{eqnarray*}
{\cal M}_{a,A}^-\left(\frac{x_N}{r'}\left(-\frac{1}{{r'}^2}X'\otimes X' +I'\right)+ \frac{x_N}{r'}(2-2\frac{x_N^2}{r^2})e_N\otimes e_N\right)&\geq& -A \frac{x_N^-}{ r^\prime} (N-1)\\
&\geq& -A (N-1)(\cot \psi)^- 
\end{eqnarray*}
and  using  $\frac{|2x_Nr'|}{r^2}\leq 1$ 
$${\cal M}_{a,A}^-(-2\frac{x_N}{r'}\frac{1}{r^2}X'\otimes X') \geq -2 {A|x_N|r^\prime\over r^2} \geq -A. $$ 
From this one gets that 

\begin{eqnarray*}
 {\cal M}^-_{a,A}(r^2D^2\theta)&\geq& -A \left((N-1)(\cot \psi)^-+1\right) \\
&& + {\cal M}^-_{a,A}(\frac{1}{r'}(1-2\frac{x_N^2}{r^2})(X'\otimes e_N+e_N\otimes X' )\\
  &\geq& -|1-2\frac{x_N^2}{r^2}|A-A \left((N-1)(\cot \psi)^-+1\right)\\
  &\geq& -A-A \left((N-1)(\cot \psi)^-+1\right)\\
  &\geq & -A \left((N-1)(\cot \psi)^-+2\right)
 \end{eqnarray*}
 where we have used that $|1-2{x_N^2\over r^2}| \leq 1$.

Putting everything together we have obtained:

\begin{eqnarray*}
{\cal M}^+_{a,A}(D^2v)&\leq& r^{\gamma-2}\left(\varphi\gamma(A(N-1)+a(\gamma-1))\right.\\
&&\left.-\varphi^\prime\left( A(N-1)(\cot \psi)^-+2)+ \gamma (A-a)\right)+a\varphi^{\prime\prime}\right)\\
&\leq & r^{\gamma-2}\left(\varphi\gamma(A(N-1)+a(\gamma-1))\right.\\
&&\left.-\varphi^\prime\left(A ((N-1)(\cot \psi)^-+2)+ \gamma (A-a)\right)+a\varphi^{\prime\prime}\right).
\end{eqnarray*}
With $\beta=A ((N-1)(\cot \psi)^-+2)+ \gamma (A-a)$ we shall choose $\varphi$ such that

$$a\varphi^{\prime\prime}-\beta\varphi^\prime+\varphi\gamma A(N-1)=0$$
and such that for $\theta$ in some interval $[0,\psi]$:

$$\varphi>0,\quad \varphi^\prime\leq 0,\quad \varphi^{\prime\prime}\leq 0.$$
Indeed, the solutions are given by
$$\varphi = C_1 e^{\sigma_1 \theta} + C_2 e^{\sigma_2 \theta}$$
with 
$\sigma_1$ and $\sigma_2$ being the positive constants 
$\sigma_1 ={1\over 2} \left( \beta + \sqrt{ \beta^2-4 \gamma \left({(N-1)A\over a}\right)}\right)$


$\sigma_2 ={1\over 2} \left( \beta - \sqrt{ \beta^2-4 \gamma \left({(N-1)A\over a}\right)}\right)$
where $\gamma$ is sufficiently close to zero in order that $ \beta ^2 >4 \gamma \left({(N-1)A\over a}\right)$. In that case $\sigma_1$ and $\sigma_2$ are both positive,  one also has  $\sigma_1> \sigma_2$. We prove that for $\gamma$ small enough, one can find a solution $\varphi$ such that on $[0, \psi]$,
$\varphi\geq 1$, $\varphi^\prime \leq 0$ and $\varphi^{\prime\prime}\leq 0$.

 We choose $C_1<0$ and $C_2>0$ with 
 $$\left\{\begin{array}{l}
 C_1\sigma_1+ C_2 \sigma_2 = 0\\
 C_1 e^{\sigma_1 \psi} + C_2 e^{\sigma _2 \psi } = 1
 \end{array}
 \right.
 $$
  This  system has a solution because 
 for $\gamma$ small enough 
 \begin{eqnarray*}
 e^{(\sigma_2-\sigma_1)\psi }& \geq &e^{-\beta\psi }\geq {4\gamma(N-1)A\over a\beta^2}\\
 &\geq & 1-\sqrt{1-{4\gamma(N-1)A\over a\beta^2}}\\
 &\geq &{\beta -\sqrt{\beta^2-{4\gamma(N-1)A\over a}}\over \beta}\\
 &\geq & {\beta -\sqrt{\beta^2-{4\gamma(N-1)A\over a}}\over \beta+\sqrt{\beta^2-{4\gamma(N-1)A\over a}}}\\
 &=& {\sigma_2\over \sigma_1}.
 \end{eqnarray*}
 We now deduce from this that 
 $\varphi^\prime \leq 0$, and $\varphi^{\prime\prime} \leq 0$ on $[0, \psi ]$.
 
 Indeed 
 the assumption implies that $\varphi^\prime (0) = 0$. 
 Then 
 for $\theta>0$
 $$\varphi^\prime (\theta) = C_1 \sigma_1 e^{\sigma_1\theta} + C_2 \sigma_2 e^{\sigma_2 \theta} \leq (C_1\sigma_1 + C_2\sigma_2) e^{\sigma_1 \theta} = 0.$$
One also has 
 $$\varphi^{\prime\prime} (0) = C_1\sigma_1^2+ C_2 \sigma_2^2 = -C_2 \sigma_1\sigma_2 + C_2 \sigma_2^2 \leq 0$$
 and  for $\theta>0$
 $$\varphi^{\prime\prime}(\theta) = C_1\sigma_1^2 e^{\sigma_1 \theta} + C_2\sigma_2^2  e^{\sigma_2 \theta} \leq e^{\sigma_1 \theta} (C_1\sigma_1^2+ C_2\sigma_2 ^2).$$

 Let us note that 
 $$1\leq \varphi(\theta) \leq \varphi(0) =: C_1+ C_2 = C_2(1-{\sigma_2\over \sigma_1})$$
 and 
 $$|\varphi^\prime (\theta)| \leq |\varphi^\prime (\psi) | = C_2 \sigma_2(e^{\sigma_1\psi-\sigma_2 \psi}).$$
 
 Let $C_{\psi} =\sup ( \varphi^2 + (\varphi^\prime)^2)^{\alpha\over 2} $.
 We have obtained that 
 \begin{eqnarray*}
  F(x, \nabla v, D^2 v)&\leq & |\grad v|^\alpha {\cal M}_{a,A}^+(D^2v)\\
  &\leq& \gamma^\alpha r^{(\gamma-1)\alpha }(\varphi^2+(\varphi^\prime)^2)^{\alpha\over 2}{\cal M}_{a,A}^+(D^2v)\\
&\leq &  -a \gamma^{2+\alpha} (1-\gamma) \varphi r^{(\gamma-1)\alpha+\gamma-2} (\varphi^2+(\varphi^\prime)^2)^{\alpha\over 2} \\
  &\leq &-a \gamma^{2+\alpha} C_{\psi} r^{\gamma(\alpha+1)-\alpha -2}.
  \end{eqnarray*}

We now consider the case
 $h\neq 0$. The above computations give 
\begin{eqnarray*}
 F(x,\nabla v, D^2 v)+ h(x)\cdot \nabla v |\nabla v|^\alpha &\leq&-C_{\psi} r^{\gamma(\alpha+1)-\alpha-2} \gamma^{2+\alpha} a \varphi\\
 &&+ |h|_\infty (\gamma r^{\gamma-1})^{1+\alpha} \sup (|\varphi|^2+( \varphi^\prime)^2) ^{1+\alpha\over 2}\\
 &<&-{C_{\psi} r^{\gamma(\alpha+1)-\alpha-2} \gamma^{2+\alpha} a \over 2}\\
 &\leq&  -{C_{\psi} r_o^{\gamma(\alpha+1)-\alpha-2} \gamma^{2+\alpha} a \over 2}\\
 &:=&-b
 \end{eqnarray*}
for $r\leq r_o:=\inf(\bar r,\frac{\gamma a}{ C_{\psi}^{1\over\alpha}|h|_\infty})$ . 
This ends the proof.

\bigskip
\begin{rema}\label{remawprime}
In the same manner one can construct a local barrier by below, i.e. some continuous non positive function $w^\prime_z$  such that $w^\prime_z(z)=0$ which in the cone is a sub-solution of
$$F[w^\prime_z] + h(x)\cdot \nabla w^\prime_z|\nabla w^\prime_z|^\alpha \geq 1.$$
\end{rema}

\subsection{Global barriers and existence.}
We now construct a global barrier which will allow us to prove the following existence result.

\begin{prop}\label{propuzero}
Let $\Omega$ and $F$ satisfy the previous assumptions.
Then, there exists $u_o$ a nonnegative viscosity solution of 
\begin{equation}\label{unot}
\left\{\begin{array}{lc}
F(x,\nabla u_o, D^2 u_o) + h(x)\cdot \nabla u_o |\nabla u_o|^\alpha = -1&\mbox{in}\quad \Omega\\
u_o=0 & \mbox{on}\quad\partial\Omega
\end{array}
\right.
\end{equation}
which  is $\gamma$ H\"older continuous. 
\end{prop}

This Proposition \ref{propuzero} will be the first step in the proof of the maximum principle and the construction of the principal eigenfunction for non smooth bounded domains.

The global barrier is given in 
\begin{prop}\label{glbar}
For all $z\in \partial \Omega$, there exists  a continuous function $W_z$ on $\overline{\Omega}$, such that 
$W_z(z)=0$,  $W_z>0$ in $\Omega \setminus \{z\}$
which is a super solution of
\begin{equation}\label{eqq1}
 F(x,\nabla W_z, D^2 W_z)+ h(x)\cdot \nabla W_z |\nabla W_z|^\alpha\leq -1\quad \mbox{in}\quad \Omega.
 \end{equation}
\end{prop}
\begin{rema} In the same way, using Remark \ref{remawprime}  one can construct a continuous function $W'_z$ on $\overline{\Omega}$, such that 
$W'_z(z)=0$,  $W'_z<0$ in $\Omega \setminus \{z\}$  which is a sub- solution of
\begin{equation}\label{eqq1}
 F(x,\nabla W'_z, D^2 W'_z)+ h(x)\cdot \nabla W'_z |\nabla W'_z|^\alpha\geq 1\quad \mbox{in}\quad \Omega.
 \end{equation}
\end{rema}

{\em Proof:} We argue on the model of \cite{CKLS}.
Choose any point $y\notin \Omega$  and $r_1$ such that 
$2r_1 < d(y, \partial \Omega)$.
Let $G_1(x) = {1\over r_1^\sigma}- {1\over |x-y|^{\sigma}}$
then 
\begin{eqnarray*}
F[G_1]+ h(x)\cdot \nabla G_1 |\nabla G_1|^\alpha &\leq& \sigma ^{1+\alpha} |x-y|^{-(\sigma+1)\alpha-\sigma-2} \left( AN -(\sigma+2) a \right.
\\
&&\left.+ |h|_\infty |x-y|\right)\\
& \leq& -(r_1)^{-\sigma (\alpha+1)-\alpha -2} \sigma^{1+\alpha}{AN \over 4}
\end{eqnarray*}
as soon as 
$$\sigma+2>  \sup ({4AN \over a}, {2|h|_\infty {\rm diam} \Omega \over  a}).$$
Moreover 

$$\frac{1}{r_1^\sigma}\geq G_1(x) \geq {2^\sigma-1\over (2r_1)^\sigma}$$
on $\overline{\Omega}$. 
Defining  $G=\frac{r_o^\gamma  r_1^\sigma}{2} G_1$,
one gets that $G\leq {r_o^\gamma\over 2} $.

We denote by $w_z(x) =|z-x|^\gamma \varphi(\theta)$ some local barrier associated to the point $z\in \partial \Omega$ as constructed in the previous section. 
Let
$$V_z(x) = \min (G(x), w_z).$$
Since the infimum of two super-solution is a super solution,  $V_z $ is a super-solution of
\begin{eqnarray*}
F[V_z]& +& h(x)\cdot \nabla V_z |\grad V_z|^\alpha\\
& \leq& -\kappa^{1+\alpha}\\
& =& \sup 
\left(-{C_{\psi} r_o^{\gamma(\alpha+1)-\alpha-2} \gamma^{2+\alpha} a \over 2},
-{r_o^\gamma\over 2}(r_1)^{-\sigma (\alpha+1)-\alpha -2} \sigma^{1+\alpha}{AN \over 4} \right).
\end{eqnarray*}
Multiplying by $\kappa$ we get that $W_z=\frac{V_z}{\kappa}$ will denote the super-solution  of (\ref{eqq1}).

\bigskip
\begin{rema} \label{remaglob} Observe that, since $G>0$ in $\overline\Omega$ there exists $\delta$ such that
$W_z(x)=\frac{w_z (x)}{\kappa}$ if $|x-z|<\delta$. Furthermore, 
by the uniform exterior cone condition  there exists $C_w>0$  such that if $|x-z|<\delta$
$$W_z (x)\leq C_w|x-z|^\gamma,$$
where $C_w$ depends on $\gamma$, $r_o$ and $\psi$ and is independent of $z\in\partial\Omega$.
\end{rema}

In the next proposition we shall see that  existence of global barriers allows to prove H\"older's regularity for solutions in non smooth domains: 
\begin{prop} \label{hold} Let
 $H_j$ be a sequence of bounded open regular sets such that 
$H_j\subset \overline{H_j}\subset H_{j+1}$, $j\geq 1$,  whose  union  equals to $\Omega$. 
Let $u_j$  be  a sequence of bounded solutions of 
$$
\left\{\begin{array}{lc}
F(x,\nabla u_j, D^2 u_j) + h(x)\cdot\nabla u_j|\nabla u_j|^{\alpha} = f_j& \mbox{in}\quad H_j\\
u_j = 0& {\rm on} \ \partial H_j.
\end{array}
\right.$$
with $f_j$ uniformly bounded. 
Then there exist $C$  and $\gamma>0$ independent on $j$  such that
$$|u_j(x)-u_j(y)|\leq C|x-y|^\gamma$$
for all $x,y\in\Omega$, where $\gamma\in (0,1)$ is given in the previous construction.
\end{prop}
{\em Proof:}  
Since $\partial H_j$ is ${\cal C}^2$, it satisfies the exterior sphere condition and a fortiori 
the exterior cone condition.   Since the $H_j$ converge to $\Omega$ which satisfies the exterior cone condition, we can choose exterior cones with opening $\psi$ and height $r_o$ which do not depend on $j$.

Using the global barriers of  Proposition \ref{glbar} and the comparison principle in $H_j$ one easily has  that for any $z\in \partial H_j$
 $$u_j\leq {|f_j|^{\frac{1}{1+\alpha}}_\infty W_z}, \quad \mbox{in}\quad H_j.$$
Let 

$$\Delta_\delta=\{(x,y) \in H_j^2\quad \mbox{such that}\quad |x-y|\leq \delta\}.$$
Let $C=\max\{\frac{2|u|_\infty}{\delta^\gamma},C_w| f_j|^{\frac{1}{1+\alpha}}_\infty\}$, we want to prove that for $\delta$ small enough, and for any $(x,y)\in \Delta_\delta$ 

\begin{equation}\label{eqhold}
u_j(x)-u_j(y)\leq C|x-y|^\gamma.
\end{equation}

In the first step we prove it on the boundary of $\Delta_\delta$.
Indeed if $|x-y|=\delta$  it is immediate from the definition of $C$. Suppose hence that $x\in H_j$ and $y\in\partial H_j$, with $|x-y|\leq \delta$. 
Then, using Remark \ref{remaglob}, for $\delta$ sufficiently small
$$u_j(x)\leq  {| f_j|_\infty^{\frac{1}{1+\alpha}} W_y}\leq C_w| f_j|_\infty^{\frac{1}{1+\alpha}}|x-y|^\gamma.$$

The second step  is to check  that the inequality (\ref{eqhold}) holds inside $\Delta_\delta$. It proceeds exactly as in the smooth case (see \cite{IL,BD2}) using hypothesis (H2) and (H3).

{\em
Proof of Proposition \ref{propuzero}.}

Let $H_j$ be a sequence of bounded open regular sets such that 
$H_j\subset \overline{H_j}\subset H_{j+1}$, $j\geq 1$,  with the union  equals to $\Omega$. 

Let $u_j$ for $j\geq 1$  be the  solution of 

$$
\left\{\begin{array}{lc}
F(x,\nabla u_j, D^2 u_j) + h(x)\cdot\nabla u_j|\nabla u_j|^{\alpha} = -1& \mbox{in}\quad H_j\\
u_j = 0& {\rm on} \ \partial H_j.
\end{array}
\right.$$
 Using the global barriers of  Proposition \ref{glbar} and the comparison principle in $H_j$ one easily has  that 
 $$u_j\leq {W_z} \quad\mbox{in}\quad H_j.$$
As a consequence,  $(u_j)_{j\geq 1}$ is  a bounded and increasing sequence - in the sense that 
$u_j\geq u_{j-1}$ on $H_{j-1}$-. 
Using Proposition \ref{hold},  the sequence  $(u_j)_j$ is  uniformly $\gamma$-H\"older continuous. 
 As a consequence  on any compact set $J\subset\Omega$,  one gets that $(u_j)_j$ converges uniformly  to some $u_o$ which satisfies 
 $$F(x,\nabla u_o, D^2 u_o)+ h(x)\cdot\nabla u_o|\nabla u_o|^\alpha= -1.$$
Furthermore $u_o$ equals $0$ on the boundary since, by passing to the limit in the previous inequality
 $$u_o\leq {W_z},$$
 for all $z\in \partial \Omega$.  We have also obtained that $u_o$ is $\gamma$ H\"older continuous.

\begin{rema} In the same manner it is possible to prove that there exists $u_o'$ a non positive $\gamma$-H\"older continuous solution  of
$$\left\{\begin{array}{lc}
F(x,\nabla u_o', D^2 u_o') + h(x)\cdot \nabla u_o' |\nabla u_o'|^\alpha = 1&\mbox{in}\quad \Omega\\
u_o'=0 & \mbox{on}\quad\partial\Omega,
\end{array}
\right.
$$
with $u_o'\geq W_z'$ for all $z\in\partial\Omega$.
\end{rema}
 
 \begin{cor}\label{corz}
Given $f\in {\cal C} (\overline{\Omega})$ there exists $u$ a $\gamma$-H\"older continuous viscosity solution of 
\begin{equation}\label{eqcorz}
\left\{\begin{array}{lc}
F(x,\nabla u, D^2 u)+ h(x) \cdot \nabla u |\nabla u|^\alpha  = f&\mbox{in}\quad\Omega,\\
u=0 &\mbox{on}\quad\partial\Omega
\end{array}
\right.
\end{equation}
with 
$$|u(x)|\leq |f|^{1\over 1+\alpha}_\infty \sup(u_o (x),-u_o'(x)). $$
Furthermore if $f\leq 0$, $u\geq 0$.
\end{cor}
{\em Proof }
Let $z_j$ be a sequence  of solutions on $H_j$ of 
$$\left\{
\begin{array}{lc}
F(x,\nabla z_j, D^2z_j)+  h(x) \cdot \nabla z_j |\nabla z_j|^\alpha = f &\mbox{in}\quad H_j,\\
 z_j = 0 &\mbox{on}\quad  \partial H_j.
 \end{array}
 \right.
 $$
By the comparison principle  on $H_j$, $u_o^\prime |f|^{1\over 1+\alpha}_\infty\leq z_j\leq  u_o|f|^{1\over 1+\alpha}_\infty $. 
Using Proposition  \ref{hold}, the sequence $(z_j)$ is uniformly $\gamma$ H\"older continuous and then $z_j$ converges on every compact set in $\Omega$ to a solution $z$ which is  $\gamma$ H\"older continuous. 

If $f\leq 0$, each $z_j$ is non-negative,  which implies that $z\geq 0$. 
Using the inequality 
$$|z_j|_\infty\leq |f |_\infty ^{1\over 1+\alpha}u_o$$
in $H_j$, one gets the final inequality by passing to the limit . 
\begin{rema} Observe that the existence of $u_o$ and $z$ solutions of (\ref{unot}) and (\ref{eqcorz}) can be done via Perron's method adapted to viscosity solutions. In particular choosing
$$u= \sup \{v, \ \mbox{subsolution of (\ref{eqcorz}) satisfying}, |f|^{1\over 1+\alpha}_\infty (W')^\star\leq v\leq |f|_\infty^{1\over 1+\alpha}W_\star\}$$
where $W_\star$ is the lower semi-continuous envelope of $\inf_{z\in\partial\Omega}W_z$ and
$(W')^\star$ is the upper semi continuous envelope of $\sup_{z\in\partial\Omega}W^\prime_z$ (The definition of viscosity solution is then  intended in the sense of semi-continuous viscosity solutions, see \cite{BD2}).  
\end{rema}
\begin{rema}\label{310}
For $V(x)\leq 0$, $u_o$ is a supersolution of 
 $$F(x,\nabla u_o, D^2 u_o)+ h(x)\cdot\nabla u_o|\nabla u_o|^\alpha+V(x)u_o|u_o|^\alpha= -1.$$
This implies that for any $f\leq 0$ there exists $u$ solution of
\begin{equation}\label{eqcorz}
\left\{\begin{array}{lc}
F(x,\nabla u, D^2 u)+ h(x) \cdot \nabla u |\nabla u|^\alpha +V(x)u|u|^\alpha = f&\mbox{in}\quad\Omega,\\
u=0 &\mbox{on}\quad\partial\Omega
\end{array}
\right.
\end{equation}
and 
$$u(x)\leq |f|^{1\over 1+\alpha}_\infty u_o (x).$$
\end{rema}
\begin{prop}\label{propK}
The function $u_o$  in Proposition  \ref{propuzero} satisfies also :  
$\forall \delta$,  there exists $K$,  a compact set in $\Omega$ such that 
$$\sup_{\Omega\setminus \overline{K}} |u_o| \leq \delta.$$
\end{prop}
{\em Proof.}
For each $z\in \partial \Omega$  we know that 
$$u_o\leq { W_z},\quad \mbox{in}\quad \Omega.$$
Let $\delta>0$ then for all $z\in \partial \Omega$  there exists $r_z$ such that for $x \in B(z, r_z)\cap\Omega$
$${W_z(x)}\leq \delta.$$
Since $\partial \Omega$ is compact one can extract from $\cup B(z,r_z)$ a finite recovering, say $\cup _{i\leq k}B(z_i, r_{z_i})$. Let then $K$ be a compact set such that 
$$\Omega \setminus K \subset \cup_{i\leq k} B(z_i , r_{z_i}).$$
We have 
$$ u_o \leq W= 
 \inf_{z_i, i\leq k} {W_{z_i}},$$
 and then 
 $u_o\leq \delta$
 in $\Omega\setminus K$.  
This ends the proof.

\begin{cor}\label{corK}
$\forall M >0$, there exists $K$ compact subset of  $\Omega$,  large enough, such that 
$$\bar\lambda (\Omega \setminus K) > M.$$
\end{cor}
{\em Proof. }
Let $\delta$ be such that $\left({1\over \delta}\right)^{1+\alpha} \geq M+|V|_\infty $,  and let $K$ be large enough in order that 
$$\sup _{\Omega \setminus K} |u_o|\leq \delta. $$
Then 
$$F[u_o] + h(x) \cdot \nabla u_o |\nabla u_o|^\alpha + (M+ V(x))u_o^{1+\alpha} = -1+ (M+ V(x)) u_o^{1+\alpha} \leq 0$$
in $\Omega\setminus K$, 
and since $u_o$ is positive one gets that $\bar\lambda (\Omega \setminus K) \geq  M$.
\subsection{Maximum principle }
\begin{defi}
We shall say that $\limsup_{x\rightarrow \partial \Omega } w(x) \leq 0$
if for all $\epsilon >0$ there exists $K$ compact in $\Omega$,  large enough in order that 
$\sup_{\Omega \setminus K} w\leq \epsilon$
\end{defi}

\begin{prop}\label{maximumprinciple}  Let $\beta(x,\cdot) $ be a nondecreasing  continuous  function such that $\beta (x,0) = 0$. 
Suppose that $w$ is uppersemicontinuous and bounded by above and satisfies 
$$F(x, \nabla w, D^2 w)+ h(x)\cdot \nabla w|\nabla w|^\alpha -\beta (x,w) \geq 0$$
with 
$$\limsup w(x_j)\leq 0$$
for all $x_j\rightarrow \partial \Omega$.
Then 
$w\leq 0$ in $\Omega$.
\end{prop}
\begin{rema} If $\beta$ is increasing then the result holds without requiring any regularity on the bounded domain $\Omega$. In that case one can use comparison principle  in Proposition \ref{comparisonprinciple}.
\end{rema}
{\em Proof:}

We assume by contradiction that $w>0$ somewhere in $\Omega$. Let $\bar x$ be a point in $\Omega $  such that $w(\bar x)>0$, and let $\gamma >0$ be such that  $\gamma u_o(\bar x) < w(\bar x)$. The function 
$$w-\gamma u_o$$ is uppersemicontinuous, bounded by above and it admits a supremum $>0$,  achieved  inside $\Omega$. Indeed, let $\epsilon < {w(\bar x)-\gamma u_o(\bar x)\over 2}$.
Let  $K$ be compact and  large enough,  in order  that $\bar x\in K$ and such that  
$w(x) \leq \epsilon$ in $\Omega \setminus K$. Then 
$(w-\gamma u_o)(x) \leq \epsilon$ in $\Omega \setminus K$. As a consequence $w-\gamma u_o$ achieves its maximum inside $K$.  
The end of the proof is the same as in the case of regular sets :  

We introduce $\psi_j(x, y) = w(x)-\gamma u_o(y) -{j\over q} |x-y|^q$. 
One can  prove  as in \cite{BD1}, that for $j$ large enough,   $\psi_j$ achieves its maximum on $(x_j, y_j)$ inside $\Omega\times \Omega $,  (more precisely  in $K\times K$), and that there exists ($X_j$, $Y_j$)  in ${\cal S}^2$ such that 
$$(j|x_j-y_j|^{q-2} (x_j-y_j), X_j)\in J^{2,+ } w(x_j)$$
$$(j|x_j-y_j|^{q-2} (x_j-y_j), -Y_j)\in J^{2,- } \gamma u_o(y_j)$$
Moreover one can choose $x_j\neq y_j$ for $j$ large enough, as it is done in \cite{BD1}. 

One has then   using (H2), (H4)  and the decreasing properties of $\beta$, 
\begin{eqnarray*}
0& \leq& F(x_j, j|x_j-y_j|^{q-2} (x_j-y_j), X_j)+ h(x_j)\cdot |x_j-y_j|^{(q-1)(\alpha+1)-1}(x_j-y_j)\\
&& -\beta (x_j, w(x_j))\\
&\leq &  F(y_j, j|x_j-y_j|^{q-2} (x_j-y_j), -Y_j)+h(y_j)\cdot |x_j-y_j|^{(q-1)(\alpha+1)-1}(x_j-y_j) \\
&&+ o(1)\\
&\leq & -\gamma^{1+\alpha}+o(1)
\end{eqnarray*}
a contradiction since $\gamma >0$.

\section{Existence of an eigenfunction}

We recall that  $V$  is some bounded and continuous function and that  $\bar\lambda(\Omega)$  is defined as :
$$\bar\lambda (\Omega)= \sup \{ \lambda, \exists \ \varphi>0\ {\rm in} \ \Omega , F[\varphi]+h(x)\cdot\nabla \varphi|\nabla \varphi|^\alpha+  (V(x)+ \lambda) \varphi^{1+\alpha} \leq 0\}.$$

\begin{theo}

Let $\Omega$ be a bounded domain which satisfies the uniform exterior cone  condition, $F$ satisfies condition (H1) to (H4) and $h$ satisfies (H5).  There exists a positive function $\phi$ solution of
$$
\left\{\begin{array}{lc}
F(x,\nabla \phi, D^2\phi) + h(x)\cdot \nabla \phi |\nabla \phi|^\alpha + (V(x)+ \overline\lambda(\Omega))\phi^{1+\alpha}=0 & \mbox{in}\quad\Omega\\
\phi=0 & \mbox{on}\quad\partial\Omega,
\end{array} 
\right.$$
which is $\gamma$-H\"older continuous. 
\end{theo}
{\em Proof }:

 Let $H_j$ be a sequence of regular subsets of $\Omega$, strictly increasing,  with union $\Omega$. 
One has for $\mu_j = \bar\lambda(H_j)$ the existence of $\phi_j>0$ an eigenfunction in $H_j$.  One  also assume that $\sup \phi_j=1$. 
Let $\mu = \lim \mu_j\geq \bar\lambda(\Omega)$.  (Let us note that the sequence $(\mu_j)$ is decreasing).

Since the $\phi_j$ are uniformly bounded, we can apply  Proposition \ref{hold} with 
$f_j= (V(x)+ \mu_j)\phi_j^{1+\alpha}$ and  we obtain that the sequence $(\phi_j)_j$ is uniformly H\"older continuous.  Up to a subsequence, the  sequence $(\phi_j)$ converges to $\phi$ a non-negative solution of
$$F(x,\nabla \phi, D^2\phi) + h(x)\cdot \nabla \phi |\nabla \phi|^\alpha + (V(x)+ \mu)\phi^{1+\alpha}=0.$$
We have to prove that $\phi$ is not identically zero.

Let  $K_1$ be  a compact set  of $\Omega $,  such that  $\bar \lambda (\Omega \setminus K_1)>\mu_1 = \bar\lambda (H_1)> \bar\lambda (\Omega)$, this is possible according to Corollary \ref{corK}. 

Let $\delta$ be small enough in order that 
$$(\bar\lambda (\Omega\setminus K_1 )+ |V|_\infty)  \delta^{1+\alpha} < 1.$$
 According to Proposition \ref{propK}, there exists $K_2$, a compact regular set, such  that  $K_1\subset K_2$   and $\sup_{\Omega \setminus K_2} u_o < \delta$.

 One has also $\bar\lambda (H_j\setminus K_2) \geq \bar\lambda (\Omega\setminus K_1) > \bar\lambda (\Omega)$.

 We observe that $u_o$ satisfies in $\Omega \setminus K_2$, (hence also in $H_j\setminus K_2$) : 
 $$F(x, \nabla u_o, D^2 u_o) +h(x)\cdot \nabla u_o|\nabla u_o|^\alpha +(\bar\lambda (\Omega\setminus K_1 ) +|V|_\infty)u_o^{1+\alpha} \leq 0$$
 which implies in particular  that
  $$F(x, \nabla u_o, D^2 u_o) +h(x)\cdot \nabla u_o|\nabla u_o|^\alpha +(\bar\lambda (H_j) + V)u_o^{1+\alpha} \leq 0.$$
On $\partial( H_j\setminus K_2)$ 
$$\phi_j\leq  1\leq {1\over \inf_{K_2} u_o}u_o,$$ 
hence  using the comparison principle Theorem \ref{thcomp}  on the  set  $H_j\setminus K_2$,  
one gets that 
\begin{equation}\label{phiuo}\phi_j\leq   \frac{1}{\inf_{K_2} u_o}u_o
\end{equation}
on $H_j\setminus K_2$. 

Let $K_3$ which contains $K_2$ such that in $\Omega\setminus K_3$,  $u_o\leq \frac{\inf_{K_2}u_o}{2}$, then 
$$\phi_j\leq   \frac{1}{2}\quad \mbox{in}\quad \Omega\setminus K_3.$$
This implies that $\sup_{K_3} \phi_j=1$,  and hence $\sup_{K_3}\phi=1$.

In particular we have obtained that $\phi$ is not zero and by the  strict maximum principle $\phi>0$ in $\Omega$. Furthermore $\mu\leq\bar\lambda(\Omega)$ and hence $\mu=\bar\lambda(\Omega)$.

Passing to the limit in (\ref{phiuo}) we also get that $\phi$ is zero on the boundary.
\section{Other   maximum
principle and eigenvalues}
In all the results of this section we still assume that $\Omega$ has the uniform exterior cone condition and $F$ and h satisfy (H1) to (H5).

We define 
$$\lambda_e = \sup\{\bar\lambda (\Omega^\prime), \ \Omega\subset \subset \Omega^\prime,\ \Omega^\prime \ {\rm regular\ and \ bounded}\}$$
and 
$$\tilde \lambda= \sup \{ \lambda, \exists \ \varphi>0\ {\rm in} \ \overline{\Omega} , F[\varphi]+h(x)\cdot\nabla \varphi|\nabla \varphi|^\alpha+  (V(x)+ \lambda) \varphi^{1+\alpha} \leq 0\}.$$

Let us note that the existence of $u_o$ implies that $\tilde \lambda >0$. 
On the other hand $\lambda_e\geq \bar\lambda>0$. 

In this section we are going to prove that 
$\lambda_e = \tilde \lambda$ 
and  that it is an "eigenvalue" in the sense that there exists some $\phi_e >0$, which satisfies 
$$
\left\{\begin{array}{lc}
F(x,\nabla \phi_e, D^2\phi_e) + h(x)\cdot \nabla \phi_e |\nabla \phi_e|^\alpha + (V(x)+ \lambda_e)\phi_e^{1+\alpha}=0 & \mbox{in}\quad\Omega\\
\phi_e=0 & \mbox{on}\quad\partial\Omega.
\end{array} 
\right.$$
Observe that clearly $\lambda_e \leq \overline{ \lambda}$,  and furthermore if $\Omega$ is smooth, the equality holds. The case where $\Omega$ is non smooth is open and the  identity of these values is equivalent to the existence of a maximum principle for $\lambda < \bar\lambda$. 
\bigskip

Let us start with the following maximum principle 
\begin{prop}
For $\lambda < \tilde \lambda $, if $w$ is a sub solution of
$$F(x, \nabla w, D^2 w)+ h(x)\cdot \nabla w|\nabla w|^\alpha+ (V(x)+ \lambda)w^{1+\alpha}\geq 0$$
satisfying
$$w(x)\leq 0\quad \mbox{on}\quad\partial\Omega$$
then $w\leq 0$ in $\Omega$.
\end{prop}
Sketch of the proof : 
Let $\varphi>0$ on $\bar\Omega$, such that 
$$ F[\varphi]+h(x)\cdot\nabla \varphi|\nabla \varphi|^\alpha+  (V(x)+ \lambda) \varphi^{1+\alpha} \leq 0.$$

\bigskip
Suppose that $w>0$ somewhere,  since $\varphi>0$ on  $\overline{\Omega}$ one can 
 define $\gamma^\prime =\sup_{x\in {\Omega}}{w\over \varphi}$ and follow the proof of \cite{BD3} to derive a contradiction. 
 
 \begin{rema}
 Observe that for $\lambda<\bar\lambda$, we don't know if the maximum principle holds when $\Omega$ is not smooth because for supersolutions satisfying 
 $\varphi>0$ in $\Omega$ we don't know if $\sup_{x\in {\Omega}}{w\over \varphi}$ is bounded.
 \end{rema}

\begin{prop}\label{phie}
 There exists $\phi_e>0$ which satisfies
 $$
\left\{\begin{array}{lc}
F[ \phi_e]+ h(x)\cdot \nabla \phi_e |\nabla \phi_e|^\alpha + (V(x)+ \lambda_e)\phi_e^{1+\alpha}=0 & \mbox{in}\quad\Omega\\
\phi_e=0 & \mbox{on}\quad\partial\Omega.
\end{array} 
\right.$$
\end{prop}
{\em Proof of Proposition \ref{phie}. }

Let $\Omega_j$ be a decreasing sequence of regular open bounded domains which contain $\Omega$. Let $\phi_j$ be some  positive eigenfunction for $\Omega_j$ such that $|\phi_j|_\infty= 1$, which exists   according to the results in \cite{BD3}. 

Using the comparison principle Proposition \ref{comparisonprinciple}, one has that for all $z\in\partial\Omega_j$
$$\phi_j\leq (|V|_\infty+\bar\lambda(\Omega_j))W^j_z$$
where $W^j_z$ is a global barrier for $\Omega_j$.
As in Remark \ref{remaglob}, $W^j_z$ satisfies 

$$W^j_z(x)\leq C|x-z|^\gamma$$
with $C$ independent of $j$ and $z$,  since the $\Omega_j$ converge to $\Omega$ which satisfies  the uniform exterior sphere condition. 
This implies that  for $\epsilon >0$ there exists $K$  compact in $\Omega$, large enough
in order that 
$$\sup_j\sup_{\Omega_j\setminus K} \phi_j\leq \epsilon.$$
 In particular $\phi_j$ has the property that if $d(K, \partial \Omega_j) < \left({\epsilon\over C}\right)^{1\over \gamma}.$
 $$\sup_{\Omega_j\setminus K} \phi_j(x)\leq \epsilon$$
 Let $K$ be a compact set in $\Omega$  such that $d(K, \partial \Omega)< \epsilon$. Since the distance is continuous, for $j$ large enough $d(K, \partial \Omega_j)< \epsilon$ and then 
 $$\sup_{\Omega_j\setminus K} \phi_j(x)\leq \epsilon.$$
 In particular one can take a compact $K$ large enough in $\Omega$ in order that 
$$\sup_{\Omega_j\setminus K} \phi_j \leq {1\over 2}$$
and then the supremum of $\phi_j$ is achieved in $K$. 

 By the uniform estimates in Proposition \ref{hold} the sequence $(\phi_j)_j$  is uniformly $\gamma$-H\"older  on $K$ and one can then extract from $(\phi_j)_j$ a subsequence such that  $\phi_j$ converges to some function $\phi_e$ which is  such that $|\phi_e|_{L^\infty (K)} = 1$. 
By compacity one has that $\phi_e$ is a solution of 
$$
F(x,\nabla  \phi_e, D^2\phi_e) + h(x)\cdot \nabla \phi_e |\nabla \phi_e|^\alpha + (V(x)+ \lambda_e(\Omega))\phi_e^{1+\alpha}=0 \  \mbox{in}\quad\Omega.$$
Moreover $\phi_e >0$ in $\Omega$, and the  estimate 
$$\phi_j \leq C \inf_{z\in \partial \Omega_j} W^j_{z}$$
gives, by passing to the limit, that $\phi_e = 0$ on the boundary of $\Omega$.

\begin{cor}

$$\lambda_e= \tilde{\lambda}$$
\end{cor}
{\em Proof:}
Suppose  by contradiction that $\lambda_e < \tilde {\lambda}$, then by the maximum principle one would obtain that $\phi_e\leq 0$. 

\bigskip
We now present some existence result for the Dirichlet problem. 

 \begin{prop} Let $\lambda<\lambda_e$ then for any function $f\leq 0$ and continuous there exists $u$  a viscosity solution of 
  \begin{equation}\label{diric}
\left\{\begin{array}{lc}
F(x,\nabla u , D^2u) + h(x)\cdot \nabla u|\nabla u|^\alpha + (V(x)+ \lambda)u^{1+\alpha}=f & \mbox{in}\quad\Omega\\
u=0 & \mbox{on}\quad\partial\Omega.
\end{array} 
\right.
\end{equation}
Furthermore $u\geq 0$ and $\gamma$-H\"older continuous.
\end{prop}
{\em Proof.} 
For $K=2|V|_\infty+|\lambda|$ let $u_n$ be the sequence of solutions of 

 $$
\left\{\begin{array}{lc}
F[u_{n+1}]+ h(x)\cdot \nabla u_{n+1} |\nabla u_{n+1} |^\alpha + (V(x)+\lambda-K) u_{n+1}^{1+\alpha} =f- K u_n^{\alpha+1}& \mbox{in}\quad\Omega\\
u_{n+1} =0 & \mbox{on}\quad\partial\Omega
\end{array} 
\right.$$
with $u_1=0$,  $u_{n}$ exists by Remark \ref{310}.
The sequence $(u_n)_n$ is increasing by the comparison principle in Proposition \ref{comparisonprinciple}. Arguing as in \cite{BD2} one can prove that the sequence is bounded, using the maximum principle of Proposition \ref{maximumprinciple}.   Furthermore there exists a constant $C$ such that
$$u_n\leq Cu_o.$$
Passing to the limit, which we can do thanks to the  H\"older's regularity given in Propostion \ref{hold}, we get the required solution.  

\begin{rema} The validity of the  maximum principle for $\lambda<\overline\lambda(\Omega)$  is  equivalent to $\lambda_e=\overline\lambda(\Omega)$ and to the existence of a solution for the Dirichlet problem (\ref{diric})  for any $\lambda<\overline\lambda(\Omega)$.
\end{rema}

\end{document}